\documentclass[a4paper,12pt,leqno]{article}
\usepackage{latexsym}
\usepackage[all]{xy}

\usepackage{amssymb} 
\usepackage{amsmath} 
\usepackage{theorem}

\def\P{{\mathbf{P}}}
\def\Z{{\mathbb{Z}}}
\def\K{{\mathbb{K}}}
\def\CC{{\mathbb{C}}}
\def\R{{\mathbb{R}}}
\def\A{{\mathcal{A}}}
\def\B{{\mathcal{B}}}

\DeclareMathOperator{\codim}{codim}

\DeclareMathOperator{\Der}{Der}

\DeclareMathOperator{\pd}{pd}

\numberwithin{equation}{section}

\newcommand{\owari}{\hfill$\square$}

\theoremstyle{break}
\newtheorem{theorem}{Theorem}[section]
\newtheorem{prop}[theorem]{Proposition}
\newtheorem{cor}[theorem]{Corollary}

\newtheorem{define}[theorem]{Definition}

\newtheorem{rem}[theorem]{Remark}

\newtheorem{example}[theorem]{Example}
\newtheorem{problem}[theorem]{Problem}

\title{Deletion theorem and combinatorics of hyperplane arrangements}
\author{Takuro Abe\footnote
{
Institute of Mathematics for Industry,
Kyushu University,
Fukuoka 819-0395, Japan.
Email:abe@imi.kyushu-u.ac.jp. 
}
}
\date{\today}

\begin{document}

\maketitle

\begin{abstract}
We show that the deletion theorem of a free arrangement is combinatorial, i.e., 
whether we can delete a hyperplane from a free arrangement 
keeping freeness depends only on the intersection lattice. 
In fact, we give an explicit sufficient and necessary condition 
for the deletion theorem in terms of characteristic polynomials. 
This gives a lot of corollaries including the 
existence of free filtrations. The proof is based on 
the result about the form of minimal generators of a logarithmic 
derivation module of a multiarrangement which 
satisfies the $b_2$-equality. 
\end{abstract}

\section{Introduction}
Let $\A$ be a central arrangement of hyperplanes in 
$V=\K^\ell$ for an arbitrary field $\K$. In this section, we use the notation in \S 2 to explain 
the background of this article, and to 
state the main results. 
In the study of hyperplane arrangements, the most 
important problem is to determine whether some property of $\A$ depends only on 
its combinatorial data (i.e., the intersection lattice $L(\A)$) or not. 
For example, when $\K=\CC$, 
the cohomology ring of the complement of $\A$ is known to be combinatorial by 
Orlik-Solomon in \cite{OS}, but the fundamental group of it is known to be not 
combinatorial by Rybnikov in \cite{Ry}. 
On the other hand, the freeness of arrangements, the most important 
algebraic property of arrangements, is not yet known to be combinatorial or not when $\ell \ge 3$. This problem is called Terao's conjecture. In general, whether some property is combinatorial or not 
is not known in most cases, and they are open problems. The aim of this article is to determine 
the deletion theorem of free arrangements is combinatorial, by giving the explicit condition for it. 
To state it, let us recall Terao's addition-deletion theorem. 

\begin{theorem}[\cite{T1}, Terao's addition-deletion theorem]
Let $H \in \A$, $\A':=\A \setminus \{H\}$ and let $\A'':=\A^H:=\{L \cap H \mid 
L \in \A'\}$. Then two of the following 
imply the third:
\begin{itemize}
\item[(1)]
$\A$ is free with $\exp(\A)=(d_1,\ldots,d_{\ell-1},d_\ell)$.
\item[(2)]
$\A'$ is free with $\exp(\A')=(d_1,\ldots,d_{\ell-1},d_\ell-1)$.
\item[(3)]
$\A''$ is free with $\exp(\A'')=(d_1,\ldots,d_{\ell-1})$.
\end{itemize}
Moreover, all the three hold true if $\A$ and $\A'$ are both free.
\label{addition}
\end{theorem}

A lot of freeness have been checked and showed by using Theorem \ref{addition}. 
In Theorem \ref{addition}, if $\A$ is free, 
then to show the freeness of the deletion $\A'$, it suffices to check 
the algebraic structure of $D(\A'')$ and the inclusion between two exponents. 
In fact, we can show that no algebra, but just a combinatorics is necessary for the 
deletion theorem. We summarize it as the main theorem in this article in the 
following.

\begin{theorem}
Let $\A$ be a free arrangement and 
$H \in \A$.
Then $\A':=\A \setminus\{H\}$ is free if and only if 
$\chi(\A^H_X;t)$ divides $\chi(\A_X;t)$ for all 
$X \in L_{i}(\A^H)$ with $2 \le i \le \ell-1$.
\label{main101}
\end{theorem}

Theorem \ref{main101} follows immediately when $\ell \le 3$ by \cite{A}, so 
it can be regarded as a higher dimensional version of it.
Though the statement in Theorem \ref{main101} is purely 
combinatorial, the proof heavily depends on algebra and algebraic geometry. 
An important corollary of Theorem \ref{main101} is the following.

\begin{cor}
Assume that $\A$ is free and take $H \in \A$. Then 
the freeness of $\A \setminus \{H\}$ depends only on 
$L(\A)$.
\label{cordeletion}
\end{cor}

Also, we have the following.

\begin{cor}
Let $\ell \le 4,\ H \in \A$ and $\chi(\A^H;t)$ divides $\chi(\A;t)$. Then 
$\A$ is free if and only if $\A^H$ is free.
\label{mOT}
\end{cor}

In \cite{AT0}, a \textbf{free filtration} of an arrangement was introduced. Namely, 
$\A$ has a free filtration 
$$
\emptyset=\A_0 \subset \cdots \subset \A_n=\A
$$
if $\{\A_i\}$ 
satisfies that $|\A_i|=i,\ |\A|=n$ and $\A_i$ is free for all $i$. 
Naively, 
a free arrangement with a free filtration is the arrangement which can be constructed from 
the empty arrangement by only using Terao's addition theorem. 
Some free arrangements are known not to have any free filtration, the most famous one is 
the cone of the affine line arrangement consisting of all edges and diagonals of 
a regular pentagon, see \cite{OT} for example. We can make it clear that this property is 
combinatorial.

\begin{cor}
For a free arrangement $\A$,
whether it has a free filtration or not depends only on 
$L(\A)$.
\label{corfreefilt}
\end{cor}







To show the results above, we need the 
\textbf{$b_2$-inequality} with respect to $H \in \A$ which was 
shown in \cite{A2} and \cite{A3}:
$$
b_2(\A) \ge b_2(\A^H)+(|\A|-|\A|)|\A^H|.
$$

It is easy to check that the $b_2$-inequality becomes a equality 
if $\chi(\A^H;t) \mid \chi(\A;t)$.
This (in)equality played a key role in 
the proof of the division theorem in \cite{A2}. 
When the $b_2$-equality holds, 
the following structure theorem on 
the logarithmic derivation modules holds, which is essential for the proof of our main 
results.

\begin{theorem}
Let $(\A,m)$ be multiarrangements such that $\A\neq \emptyset$ and that 
$$
b_2(\A,m)=b_2(\A)+(|m|-|\A|)(|\A|-1).
$$
Then 
\begin{itemize}
\item[(1)]
there is a minimal generator $\theta_E,\theta_1,\ldots,\theta_s$ for 
$D(\A)$ such that 
$$
\frac{Q(\A,m)}{Q(\A)} \theta_E, \theta_1,\ldots,\theta_s
$$
form a generator for $D(\A,m)$, 
and they form a minimal generator for $D(\A,m)$ 
if and only if 
$\langle \theta_1,\ldots,\theta_s \rangle_S \neq D(\A,m)$,
\item[(2)]
if $(\A,m)$ is free and $\A$ is not free, 
then 
there exists a minimal generator $\theta_E,\theta_1,\ldots,\theta_\ell \in D(\A)$ for $D(\A)$ such that 
$\theta_1,\ldots,\theta_\ell$ form a basis for $D(\A,m)$, and 
\item[(3)]
if $(\A,m)$ and $\A$ are both free, then 
there is a free basis $\theta_E,\theta_1,\ldots,\theta_{\ell-1}$ for 
$D(\A)$ such that $\theta_1,\ldots,\theta_{\ell-1},(Q(\A,m)/Q(\A) )\theta_E$ 
form a basis for $D(\A,m)$.
\end{itemize}
\label{main2}
\end{theorem}

The most important application of Theorem \ref{main2} is the following.

\begin{cor}
Let $\A$ be an $\ell$-arrangement and 
$(\A^H,m^H)$ the Ziegler restriction of $\A$ onto $H \in \A$. Assume that 
the $b_2$-equatity
$$
b_2(\A)=b_2(\A^H)+|\A^H|(|\A|-|\A^H|)
$$
holds true. Let 
$\theta_E^H:=(Q(\A^H,m^H)/Q(\A^H))\theta_E \in D(\A^H,m^H)$, and 
let $\pi:D_H(\A) \rightarrow D(\A^H,m^H)$ be the Ziegler 
restriction.

(1)\,\,
Then there are minimal generators $\theta_E,\theta_1,\ldots,\theta_s$ for 
$D(\A^H)$ such that $\theta_E^H,
\theta_1,\ldots,\theta_s$ form a generator for $D(\A^H,m^H)$. 

(2)\,\,
Assume that $\A$ 
is free, and $\A^H$ is not free. Then there is a generator $
\theta_E,\theta_2,\ldots,\theta_\ell$ for $D(\A^H)$ such that $\theta_2,\ldots,
\theta_\ell$ form a basis for $D(\A^H,m^H)$, 
the preimage of $\theta_2,\ldots,\theta_\ell$ in $D_H(\A)$ by $\pi$ 
form a free basis for
$D_H(\A)$, and the relation among $\theta_E,\theta_2,\ldots,\theta_\ell$ in 
$D(\A^H)$ 
is in the degree 
$d:=|\A|-|\A^H|$ of the form 
$
\theta_E^H=\sum_{i=2}^\ell f_i \theta_i
$
for $f_i \in S/\alpha_H S$, 
and no other relation exists. In other words, we have a 
free resolution
$$
0 \rightarrow 
S[-d] \rightarrow 
\oplus_{i=1}^\ell S[-d_i] \rightarrow D(\A^H) \rightarrow 0,
$$
where $\deg \theta_i=:d_i$.
\label{cormain}
\end{cor}

We investigate several properties of $D(\A)$ by using Theorem \ref{main2} including 
the modified Orlik's conjecture (Problem \ref{Orlik}). 

The organization of this article is as 
follows: In \S 2 we introduce a notation and several results used for the proof. 
In \S 3 we prove Theorem \ref{main2} and show several related results. 
In \S 4 we prove main results. 
\medskip

\noindent
\textbf{Acknowledgements}. 
The author is 
partially supported by JSPS Grant-in-Aid for Scientific Research (B) JP16H03924, and 
Grant-in-Aid for Exploratory Research JP16K13744. We are grateful to Michael 
DiPasquale for this informing a counter example to Problem \ref{nf} which is now in 
Remark \ref{nfans}.

\section{Preliminaries}
In this section let us summarize several definitions and results used in this 
article. We refer \cite{OT} for a general reference in this section. 
Let $\K$ be an 
arbitrary field and $\A$ a 
\textbf{central arrangement of hyperplanes} in $V=\K^\ell$, i.e., 
a finite set of linear hyperplanes in 
$V$. 
Assume that 
every hyperplane $H \in \A$ is defined by a linear form $\alpha_H=0$. Let 
$Q(\A):=\prod_{H \in \A} \alpha_H$. Without any 
specification, we assume that $\A \neq \emptyset$. Let 
$S=\K[x_1,\ldots,x_\ell]$ be a coordinate ring of $V$ and $\Der S:=
\oplus_{i=1}^\ell S \partial_{x_i}$. Then the \textbf{logarithmic vector field} $D(\A)$ of $\A$ is defined as 
$$
D(\A):=\{\theta \in \Der S \mid \theta(\alpha_H) \in S \alpha_H\ (\forall H \in \A)\}.
$$
$D(\A)$ is a reflexive module, and not free in general. We say that $\A$ is \textbf{free} with $\exp(\A)=(d_1,\ldots,d_\ell)$ if 
$D(\A)$ is a free $S$-module of rank $\ell$ with homogeneous basis $\theta_1,\ldots,\theta_\ell$ such that 
$\deg \theta_i=d_i$. Since $\A$ is not empty, the lowest degree basis element $\theta_1$ can be chosen as the \textbf{Euler 
derivation} $\theta_E=\sum_{i=1}^\ell x_i \partial_{x_i}$ which is always contained in $D(\A)$. Also, for $H \in \A$, define 
$D_H(\A):=\{\theta \in D(\A) \mid \theta(\alpha_H)=0\}$.

Next let us introduce combinatorics and topology of arrangements. 
Let 
$$
L(\A):=\{ \cap_{H \in \B} H \mid \B  \subset \A\}
$$ 
be the 
\textbf{intersection lattice} of $\A$ with a partial order induced from the 
reverse inclusion. Define 
$L_i(\A):=
\{X \in L(\A) \mid 
\codim_V X=i \}$. 
The \textbf{M\"{o}bius function} $\mu : L(\A) 
\rightarrow \Z$ is defined by $\mu(V)=1$, and by 
$\mu(X):=-\sum_{X  \subsetneq Y \subset V} \mu(Y)$ for $L(\A) \ni X \subsetneq V$. 
Define the \textbf{characteristic polynomial} $\chi(\A;t)$ by 
$$
\chi(\A;t):=\sum_{X \in L(\A)} \mu(X) t^{\dim X},
$$
and the \textbf{Poincar\'e polynomial} $\pi(\A;t)$ by 
$$
\pi(\A;t):=\sum_{X \in L(\A)} \mu(X) (-t)^{\codim X}.
$$
For $X \in L(\A)$, the \textbf{localization} 
$\A_X$ of $\A$ at $X$ is defined by 
$$
\A_X:=\{H \in \A \mid 
X \subset H\},
$$
and 
the 
\textbf{restriction} 
$\A^X$ of $\A$ onto $X$ is defined by 
$$
\A^X:=\{H \cap X \mid
H \in 
\A \setminus \A_X\}.
$$
It is easy to check that $\A_X$ is free if $\A$ is free for any $X \in L(\A)$. 
Also, we say that $\A$ is \textbf{locally free} if 
$\A_X$ is free for any $X \in L(\A)$ with $X \neq \{0\}$. $\A$ is locally free if and only if 
the sheaf $\widetilde{D_H(\A)}$ is a vector bundle on $\P^{\ell-1} =\P(V^*)$ for 
any $H \in \A$.
Define the \textbf{Euler restriction map} $\rho:D(\A) \rightarrow 
D(\A^H)$ by taking modulo $\alpha_H$. Then it is known (e.g., see 
\cite{OT}) that there is an exact sequence 
$$
0 \rightarrow D(\A \setminus \{H\})
\stackrel{\cdot \alpha_H}{\rightarrow}
D(\A ) 
\stackrel{\rho}{\rightarrow} D(\A^H).
$$
The most useful 
inductive method to compute $\chi(\A;t)$ 
is so called the \textbf{deletion-restriction formula} as follows:
$$
\chi(\A;t)=\chi(\A';t)-\chi(\A^H;t).
$$
We may apply this to compute $\chi(\A;t)$ efficiently. 

Let $\chi(\A;t)=\sum_{i=0}^\ell (-1)^i b_{i}(\A) t^{\ell-i}$. 
When $\A \neq \emptyset$, it is known that $\chi(\A;t)$ is divisible by 
$t-1$. Define $\chi_0(\A;t):=\chi(\A;t)/(t-1)=\sum_{i=0}^{\ell-1} (-1)^i b_i^0(\A) t^{\ell-i-1}=
\sum_{i=0}^{\ell-1} (-1)^i b_i(d\A) t^{\ell-i-1}$, 
where $d\A$ is the deconing of $\A$ by any line $H \in \A$. It is known that 
$b_i(\A)$ is the $i$-th Betti number of $V \setminus \cup_{H \in \A} H$ when 
$\K=\CC$. Then we may relate the exponents of free arrangements and the combinatorics and 
topology as follows:

\begin{theorem}[Terao's factorization, \cite{T2}]
Assume that $\A$ is free with $\exp(\A)=(d_1,\ldots,d_\ell)$. Then 
$\chi(\A;t)=\prod_{i=1}^\ell(t-d_i)$.
\label{Teraofactorization}
\end{theorem}

A pair $(\A,m)$ is a \textbf{multiarrangement} if 
$m:\A \rightarrow \Z_{\ge 1}$. Let 
$|m|:=\sum_{H \in \A} m(H)$ and $Q(\A,m):=
\prod_{H \in \A} \alpha_H^{m(H)}$. For two multiplicities 
$m$ and $k$ on $\A$, $ k \le m$ means that 
$k(H) \le m(H)$ for all $H \in \A$. We may define its 
\textbf{logarithmic derivation module} 
$D(\A,m)$ as 
$$
D(\A,m):=\{\theta \in \Der S \mid \theta(\alpha_H) \in S \alpha_H^{m(H)}\ (\forall H \in \A)\}.
$$
We may define the \textbf{freeness} and \textbf{exponents} of $(\A,m)$ 
in the same way as 
for $m \equiv 1$. 
Also, we can define the characteristic polynomial $\chi(\A,m;t)=\sum_{i=0}^\ell (-1)^i b_i(\A,m)t^{\ell-i}$ of $(\A,m)$ in algebraic way, see 
\cite{ATW} for details. Now let us introduce the fundamental method to determine the 
freeness of $(\A,m)$.

\begin{theorem}[Saito's criterion, \cite{Sa}, \cite{Z}]
Let $\theta_1,\ldots,\theta_\ell$ be a homogeneous element in $D(\A,m)$. Then 
$D(\A,m)$ has a free basis $\theta_1,\ldots,\theta_\ell$ if and only if 
they are $S$-independent, and 
$|m|=\sum_{i=1}^\ell \deg \theta_i$.
\label{Saito}
\end{theorem}

We can construct the multiarrangement canonically from an arrangement $\A$ in the 
following manner:

\begin{define}[\cite{Z}]
For an arrangement $\A$ in $\K^\ell$ and $H \in \A$, define the \textbf{Ziegler 
multiplicity} $m^H:\A^H \rightarrow \Z_{>0}$ by  
$m^H(X):=|\{L \in \A \setminus \{H\} \mid 
L \cap H=X\}|$ for $X \in \A^H$. The pair $(\A^H,m^H)$ is 
called the \textbf{Ziegler restriction} of $\A$ onto $H$. Also, there is a canonical 
\textbf{Ziegler restriction map}
$$
\pi:D_H(\A) \rightarrow D(\A^H,m^H)
$$
by taking modulo $\alpha_H$.
\label{Zieglerrest}
\end{define}

A remarkable property of Ziegler restriction maps is the following.

\begin{theorem}[\cite{Z}]
Assume that $\A$ is free with $\exp(\A)=(1,d_2,\ldots,d_\ell)$. Then 
for any $H \in \A$, the Ziegler restriction $(\A^H,m^H)$ is also free with 
$\exp(\A^H,m^H)=(d_2,\ldots,d_\ell)$. In particular, $\pi$ is surjective when $\A$ is free.
\label{Ziegler2}
\end{theorem}

Moreover, a converse of Theorem \ref{Ziegler2} holds true with 
additional conditions.

\begin{theorem}[Yoshinaga's criterion, \cite{Y2}, \cite{AY}]
In the notation of Definition \ref{Zieglerrest}, it holds that 
$$
b^0_2(\A) \ge b_2(\A^H,m^H).
$$
Moreover, $\A$ is free if and only if the above inequality is the equality, and 
$(\A^H,m^H)$ is free. 
\label{Y}
\end{theorem}

An immediate consequense of Theorem \ref{Y} with the addition-deletion theorem 
for multiarrangements in \cite{ATW2} is the following inequality, which also 
induces the improvement of Terao's addition theorem, and the inequality is the key 
of this article.

\begin{theorem}[$b_2$-inequality and the division theorem, \cite{A2}]
It holds that 
$$
b_2^0(\A) \ge b_2(\A^H)+(|\A^H|-1)(|\A|-|\A^H|-1),
$$
which is equivalent to 
$$
b_2(\A) \ge b_2(\A^H)+|\A^H|(|\A|-|\A^H|).
$$
The equality holds if and only if $\A_X$ 
is free for all $X \in L(\A^H)$ with $\codim_V X=3$. This equality holds true if 
$\chi(\A^H;t) \mid \chi(\A;t)$. Moreover, $\A$ is free if the $b_2$-inequality is 
an equality, and $\A^H$ is free for some $H \in \A$.
\label{division}
\end{theorem}



\begin{define}
For $H \in \A$, define the derivation 
$\theta_E^H \in D(\A^H,m^H)$ by 
$$
\theta_E^H:=\frac{Q(\A^H,m^H)}{Q(\A^H)} \rho(\theta_E).
$$
\label{thetaEH}
\end{define}

Not only the Ziegler restriction in Definition \ref{Zieglerrest}, but 
also we have the other restriction, called the \textbf{Euler restriction}. See 
\cite{ATW2}, Definition 0.2 for details. Let $(\A^H,m^*)$ be the 
Euler restriction of $(\A,m)$ onto $H$. Then we have the following.

\begin{prop}[\cite{ATW2}, Definition 3.3, Lemma 3.4]
Let $(\A^H,m^*)$ be the Euler restriction of $(\A,m)$ onto $H \in \A$. 
Then there is a polynomial $B$ of degree $|m|-|m^*|-1$ such that 
$$
\theta(\alpha_H) \in (\alpha_H^{m(H)-1}, B).
$$
\label{B}
\end{prop}

\section{Proof of Theorem \ref{main2} and related results}

First let us show Theorem \ref{main2}, which will play the key roles in the 
rest of this article.
\medskip

\noindent
\textbf{Proof of Theorem \ref{main2}}.
(1)\,\,
Let 
$\theta_E,\theta_1,\ldots,\theta_s$ be a minimal generator for $D(\A)$. 
Let $1 \le m' < m$.
We show that 
$D(\A,m')$ has a generator 
$$
\varphi_{m'}:=\frac{Q(\A,m')}{Q(\A)} \theta_E, \theta_1',\ldots,\theta_s',
$$
where $\theta_i'=\theta_i+f_i \theta_E$ for some $f_i \in S$. 
We show by the induction on $|m'|$. When $m' \equiv 1$, then there is nothing to show.
Assume that the statement holds true for $m'$, and we show the same is true for 
$k:=m'+\delta_L$ with  
$L \in \A$ such that 
$k \le m$, where $\delta_L(K)=1$ when $K=L$, and $0$ otherwise.

By Proposition \ref{B}, there is a homogeneous polynomial 
$B$ of degree $|m'|-|m^*|$ such that, for any $\theta \in D(\A,m')$, 
$$
\theta(\alpha_L) \in (\alpha_L^{m'(L)+1},B),
$$
where $(\A^L,m^*)$ is the Euler restriction of $(\A,m')$ onto $L$.  
Note that $b_2(\A,k)=b_2(\A)+(|k|-|\A|)(|\A|-1)$ 
since we have the $b_2$-equality $b_2(\A,m)=b_2(\A)+
(|\A|-1)(|m|-|\A|)$, and $m' <k \le  m$.
Hence the proof of Theorem 1.7 in \cite{A2}, or Proposition 2.5 in \cite{A3} 
implies that 
$(\A^L,m^*)=(\A^L, m^L)$, where the latter is the Ziegler restriction of 
$\A$ onto $L$. Hence $|m^*|=|m^L|=|\A|-1$. In particular,
$$
|m'|-|m^*|=|m'|-|\A|+1=
\deg B=\deg \varphi_{m'}.
$$
Note that 
$\varphi_{m'}(\alpha_L) \not 
\in (\alpha_L^{m'(L)+1})$ by definition. Hence we may assume that 
$$
\varphi_{m'}(\alpha_L) \equiv B\ (\mbox{mod}\ \alpha_L^{m'(L)+1}).
$$
Hence when $\theta_i'(\alpha_L)=a_i \alpha_L^{m'(L)+1} +b_i B$ for 
some $a_i,b_i \in S$, replacing $\theta_i'$ by 
$\theta_i'':=\theta_i'-b_i \varphi_{m'}$, it holds that 
$\langle \varphi_{m'},\theta_1'',\ldots,\theta_s''\rangle_S=D(\A,m')$ and that 
$$
\alpha_L \varphi_{m'},\theta_1'',\ldots,\theta_s'' \in D(\A,k).
$$
Conversley, let $\theta \in D(\A,k)$. Since $D(\A,k) \subset D(\A,m')$, it holds that 
$$
\theta=g \varphi_{m'}+\sum_{i=1}^s g_i \theta_i'',
$$
here $g,g_i \in S$ and we used the fact that $\varphi_{m'},\theta_1'',\ldots,\theta_s''$ form a generator for 
$D(\A,m')$ too. Hence 
$$
\theta(\alpha_L)=g \varphi_{m'}(\alpha_L)+\sum_{i=1}^s g_i \theta_i''(\alpha_L)
$$
for some $g,g_i \in S$. Note that 
$\alpha_L^{m'(L)+1} \mid \theta(\alpha_L)$ and $
\alpha_L^{m'(L)+1} \mid \theta_i(\alpha_L)$ for $i=1,\ldots,s$. Since 
$\alpha_L^{m'(L)+1} \nmid \varphi_{m'}(\alpha_L)$ by the definition, 
it holds that $\alpha_L \mid g$, which shows that 
$\alpha_L \varphi_{m'},\theta_1'',\ldots,\theta_s''$ form a generator for $D(\A,k)$. The minimality is clear by the construction, 
which 
completes the proof.

(2)\,\,
Assume that $\A$ is not free, $(\A,m)$ free and $s > \ell$. By the above, there is a minimal generator 
$\theta_E,\theta_1,\ldots,\theta_s$ for $D(\A)$ such that 
$Q'\theta_E,\theta_1,\ldots,\theta_s$ form a generator for $D(\A,m)$, 
where $Q':=\frac{Q(\A,m)}{Q(\A)}$. Since 
$D(\A,m)$ is free with rank $\ell$, $Q'\theta_E,\theta_1,\ldots,\theta_s$ is not a 
minimal generator for $D(\A,m)$. Assume 
that $\theta_s$ is removable. Then $\theta_E,\theta_1,\ldots,\theta_{s-1}$ form a 
generator for $D(\A)$, which contradicts the minimality of the 
generator. Hence $Q'\theta_E$ is removable, and no $\theta_i$ 
is. Hence $\theta_1,\ldots,\theta_s$ has to be a minimal generator (thus a free basis)
for $D(\A,m)$, which contradicts $s > \ell$. 
Hence $s=\ell$, which completes the proof. 

(3)\,\,
Assume that $\A$ and $(\A,m)$ are both free. Then by (1), there 
is a generator (in fact, a basis) $\theta_E,\theta_1,\ldots,\theta_{\ell-1}$ for 
$D(\A)$ such that $Q'\theta_E,\theta_1,\ldots,\theta_{\ell-1}$ form a generator 
for $D(\A,m)$, which completes the proof. 
\owari
\medskip

\noindent
\textbf{Proof of Corollary \ref{cormain}}.
(1)\,\, Immediate from Theorem \ref{main2} since 
the $b_2$-equality holds, and 
$$
b_2^0(\A) \ge b_2(\A^H,m^H) \ge 
b_2(\A^H)+(|\A^H|-1)(|\A|-|\A^H|-1)
$$
by Theorems \ref{Y} and \ref{division}.

(2)\,\,
Let $\theta_E,\theta_2,\ldots,\theta_\ell$ be a minimal generator for $D(\A^H)$ 
as in Theorem \ref{main2} (3). 
Then $\theta_2,\ldots,\theta_\ell \in D(\A^H,m^H)$ form a basis for $D(\A^H,m^H)$ by Theorem \ref{main2} (3), and 
$\pi$ is surjective by Theorem \ref{Ziegler2}. 
Let $\varphi_i \in D_H(\A)$ be a preimage of $\theta_i$ by $\pi$. Assume that 
$\varphi_2,\ldots,\varphi_\ell$ are $S$-dependent, i.e., 
$$
\sum_{i=2}^\ell g_i \varphi_i=0
$$
for some $g_i \in S$. We may assume that not all $g_i$'s are divisible by $\alpha_H$. Then 
this relation contradicts the $S/\alpha_H$-independency of $\theta_2,\ldots,\theta_\ell$. 
Hence they are $S$-independent, and form a basis for $D_H(\A)$ by the reason 
of degrees and Saito's criterion. Since $\theta_E^H \in D(\A^H,m^H)$, and 
$\theta_2,\ldots,\theta_\ell$ form a basis for $D(\A^H,m^H)$, there are 
$f_i \in S/\alpha_H$ such that 
$$
\theta_E^H=\sum_{i=2}^\ell f_i \theta_i
$$
in degree $|\A|-|\A^H|$. Assume that there is $\theta \in D(\A^H)$ such that 
$\theta \mid \theta_E^H$, and 
$$
\theta=\sum_{i=2}^\ell f_i' \theta_i
$$
for some $f_i' \in S/\alpha_H$. Since the right hand side is in $D(\A^H,m^H)$, 
so is the left hand side. Thus $\theta=\theta_E^H$. Also, there are no relation
among $\theta_2,\ldots,\theta_\ell$ since they are basis for $D(\A^H,m^H)$, which 
completes the proof.\owari
\medskip

The following is immediate.

\begin{cor}
Assume that $\A$ is free and the $b_2$-equality holds for $H \in \A$. Then 
$D(\A^H)$ is generated by at most $\ell$-elements.
\label{gen}
\end{cor}

\begin{example}
Recall the Edelman-Reiner's arrangement $\A$ in $\R^5$ in \cite{ER} defined by 
$$
(\prod_{i=1}^5 x_i) (x_1\pm x_2 \pm x_3 \pm x_4 \pm x_5)=0.
$$
This consists of $21$-hyperplanes. Let $H$ be 
$x_1-x_2-x_3-x_4-x_5=0$. Then Edelman and Reiner showed in \cite{ER} that 
$\A$ is free with $\exp(\A)=(1,5,5,5,5)$ but $\chi_0(\A^H;t)=(t-4)(t-10t+26)$, hence 
not free. However, since $b_2^0(\A)=150$ and 
$b_2(\A^H)=80$, it is clear that the $b_2$-equality holds for $H \in \A$. Hence we may 
apply Corollary \ref{gen} to show that 
$\A^H$ is generated by the Euler derivations and $4$-derivations of 
degree $5$, and this is a minimal generator. In particular, since 
the relation between them are in degree $6=|\A|-|\A^H|$ by Theorem \ref{main2}, 
we can see that $\A^H$ is a nearly free surface in the sense of \cite{DiS1} and 
\cite{DiS2}. Also, since the known counter examples to Orlik's conjecture 
are very few, we may pose the following problem.
\label{ER1}
\end{example}

\begin{problem}
Let $\A$ be a free arrangement. Then is $\A^H$ either free or nearly free for $H \in \A$?
\label{nf}
\end{problem}


\begin{rem}
To Problem \ref{nf}, we learned a counter example from Michael DiPasquale as follows. 
Let $\A$ be defined by 
$$
(x_1^2-x_0^2)
(x_2^2-x_0^2)
(x_3^2-x_0^2)
(x_4^2-x_0^2)
(x_1-x_2)
(x_2-x_3)
(x_3-x_4)
(x_4+x_1)
x_0=0
$$
in $\R^5$. Let $H:=\{x_0=0\} \in \A$. Then 
$\A$ is free with $\exp(\A)=(1,3,3,3,3)$, but $\A^H$ is not free with $\pd_S D(\A^H)=2$, which 
implies that $\A^H$ is not nearly free in the sense of \cite{DiS2}. 
\label{nfans}
\end{rem}

\section{Proof of the main results}

Now let us go to the proof of Theorem \ref{main101}. The following is 
the starting point of the proof.

\begin{theorem}
Let $\A$ be a free arrangement, 
$H \in \A$ and $\chi(\A^H;t)$ divides $\chi(\A;t)$. 
If $\A^H$ is locally free, then $\A^H$ is free. In particular, by Terao's addition-deletion 
theorem, $\A \setminus \{H\}$ is free too.
\label{main1000}
\end{theorem}

\noindent
\textbf{Proof}. 
Assume that $\A^H$ is not free. Let $\exp(\A)=(d_1,\ldots,d_\ell)$ and 
$|\A|-|\A^H|=d_\ell$. 
By Corollary \ref{cormain} (2), 
$D(\A^H)$ is generated by 
$\rho(\theta_E),\pi(\theta_2),\ldots,\pi(\theta_\ell)$, where 
$\theta_2,\ldots,\theta_\ell$ form a basis for $D_H(\A)$ of 
$\deg \theta_i=d_i$, $\pi:D_H(\A) \rightarrow D(\A^H)$ the 
Ziegler restriction and $\rho:D(\A) \rightarrow D(\A^H)$ the Euler 
restriction. Since $\A$ is free, 
$D(\A^H,m^H)$ is free with basis 
$\pi(\theta_2),\ldots,\pi(\theta_\ell)$, and by Corollary \ref{cormain} (2), 
$\theta_E^H:=Q(\A^H,m^H)/Q(\A^H) \rho(\theta_E)$ 
can be expressed as a linear combination of 
$\pi(\theta_2),\ldots,\pi(\theta_\ell)$ as follows:
$$
\theta_E^H=\sum_{i=2}^\ell f_i \pi(\theta_i)
$$
for $f_i \in S/\alpha_H$. Note that 
$$
\deg \theta_E^H=|\A|-1-|\A^H|+1=d_\ell.
$$
Assume that some $f_i \neq 0$ for some $i$ with $
d_i=d_\ell$. Then it is clear that 
$$
\pi(\theta_2),\ldots,
\pi(\theta_{i-1}),\theta_E^H,\pi(\theta_{i+1}),\ldots,\pi(\theta_\ell)
$$
are $S/\alpha_H$-independent, and the sum of degrees coincides with 
$|m^H|$. Hence Saito's criterion shows that they form a basis for $D(\A^H,m^H)$. Moreover, 
the $S/\alpha_H$-linear derivations 
$$
\pi(\theta_2),\ldots,
\pi(\theta_{i-1}),\rho(\theta_E),\pi(\theta_{i+1}),\ldots,\pi(\theta_\ell)
$$
form a basis for $D(\A^H)$ since the sum of degrees coincide with 
$|m^H|-(\deg Q(\A^H,m^H)-\deg Q(\A^H))=|\A^H|$. Hence 
$\A^H$ is free, a contradiction. Assume that $f_i=0$ if 
$d_i=d_\ell$. Let $\theta_1:=\theta_E$. 
Then by Corollary \ref{cormain} (2), there is a resolution
$$
0 \rightarrow 
\mathcal{O}_H(-d_\ell) \stackrel{f}{\rightarrow}
\oplus_{i=1}^\ell \mathcal{O}_H(-d_i) \rightarrow 
\widetilde{D(\A^H)} \rightarrow 0.
$$
Here $f =(f_1,\ldots,f_\ell) \in 
\oplus_{i=1}^\ell H^0(\mathcal{O}_H(d_\ell-d_i))$. 
By the above argument, $f_i=0$ if $d_i \ge d_\ell$. Let 
$I:=\{i\mid 2 \le i \le \ell,\ d_i<d_\ell\}$. By definition, 
$|I|\le \ell-2$, thus 
$$
\P^{\ell-2}=\P(H) \supset Z_f:=\cap_{i \in I} \{f_i=0\} \neq \emptyset. 
$$ 
Recall the form of the relation:
$$
\frac{Q(\A^H,m^H)}{Q(\A^H)} 
\rho(\theta_E)=f_1 \rho(\theta_E)=
\sum_{i=2}^\ell f_i \theta_i.
$$
Since $\rho(\theta_E) \neq 0$ for any points $x \in \P(H) \simeq 
\P^{\ell-2}$, $f_2=\cdots=f_\ell=0$ at some $x \in \P(H)$ implies that 
$f_1=0$ at $x$. Hence $Z_f \cap \{f_1=0\}=Z_f\neq \emptyset $.
In other words, there is a point $x \in \P(H)$ such that 
$f_i=0$ at $x$ for all $i$. Now we show that this contradicts the assumption that 
$\A^H$ is locally free.  Take $x \in Z_f$ and consider the following 
exact sequence:
$$
\mathcal{O}_H(-d_\ell)_x \otimes k_x \stackrel{f}{\rightarrow}
\oplus_{i=1}^\ell \mathcal{O}_H(-d_i)_x \otimes k_x \rightarrow 
\widetilde{D(\A^H)}_x \otimes k_x \rightarrow 0,
$$
where $k_x$ is the residue field of the stalk $\mathcal{O}_{x,H}$. The above is 
$$
k_x \stackrel{f}{\rightarrow}
k_x^\ell \rightarrow 
k_x^{\ell-1} \rightarrow 0.
$$
By the choice of $x$, $f$ is a zero map at $x$, 
a contradiction.\owari
\medskip

\noindent
\textbf{Proof of Theorem \ref{main101}}. 
First let us show the ``if'' part by 
induction on $\ell$. By Theorem \ref{main1000}, 
this is true when $\ell \le 4$. Assume that $\ell \ge 
5$. Note that $\A_X$ is free since $\A$ is free. Then the assumption that 
$\chi(\A^H_X;t)$ divides $\chi(\A_X;t)$ for all 
$X \in L_{i}(\A^H)$ with $2 \le i \le \ell-1$ implies that 
$\A_X \setminus \{H\}$ is free by the freeness of $\A_X$ and the induction hypothesis. 
Hence $\A^H$ is 
locally free by Theorem \ref{addition}, thus Theorem \ref{main1000} completes 
the proof. 

Next let us show the ``only if'' part. 
Since $\A\setminus \{H\}$ is free, 
both $\A_X$ and $(\A \setminus \{H\})_X=\A_X \setminus \{H\}$ are free too.
Hence Theorem \ref{addition} says that $\A_X^H$ is also free, and 
$\chi(\A_X^H;t)$ divides $\chi(\A_X;t)$ by Theorem \ref{Teraofactorization}.
\owari
\medskip

Now Corollaries \ref{cordeletion} and \ref{mOT} follow immediately. 
\medskip

\noindent
\textbf{Proof of Corollary \ref{corfreefilt}}.
Induction on $|\A|$. It is trivial if $|\A|\le 1$. Assume that $|\A|>1$. 
If $\A$ has a free filtration, then at least one hyperplane can be removable 
from $\A$ keeping freeness, which is combinatorial by Theorem \ref{main101}. 
Then apply the induction hypothesis to the deleted free arrangements.\owari
\medskip

Also from the proof, the following is 
immediate too.

\begin{cor}
Let $\A$ be free and $H\in \A$. Then $\A \setminus \{H\} $ is free if 
and only if the $b_2$-equality holds for $H \in \A_X$, and 
$|\A_X|-|\A^H_X|$ is a root of $\chi(\A_X;t)$ for all $X \in L(\A^H)$.
\label{b1}
\end{cor}

Now let us relate the results above to the modified Orlik's 
conjecture. 
Orlik's conjecture asserts that the restriction 
$\A^H$ of a free arrangement $\A$ onto $H \in \A$ is free, which was 
settled negatively by Edelman and Reiner in \cite{ER}. Based on the division
theorem in \cite{A2}, the author posed the following modified Orlik's conjecture:

\begin{problem}[\cite{A3}]
Let $\A$ be a free arrangement, $H \in \A$ and 
$\chi(\A^H;t) \mid \chi(\A;t)$. Then $\A^H$ is free.
\label{Orlik}
\end{problem}

If we replace the freeness of $\A$ by that of $\A^H$ in Problem \ref{Orlik}, then 
the statement is true by Theorem \ref{division}. We cannot show Problem \ref{Orlik} without 
any assumption. What we can say from the main results is as follows.

\begin{theorem}
Problem \ref{Orlik} is combinatorial, i.e., in the assumption in Problem  
\ref{Orlik}, whether $\A^H$ is free or not depends only on $L(\A)$.
\label{combOrlik}
\end{theorem}

\begin{cor}
The modified Orlik's conjecture holds true for $\ell \le 4$.
\label{OC}
\end{cor}

The following gives a sufficient condition for the freeness of the modified Orlik's 
conjecture.

\begin{cor}
Assume that $\A$ is free with $\exp(\A)=(d_1,\ldots,d_\ell)$, $H\in \A$ and 
$\chi(\A^H;t) $ divides $\chi(\A;t)$. 
Then $\A^H$ is free if  there is a 
free subarrangement $\B \subset \A$ such that $\B \ni H$ and $\B^H=\A^H$.
\label{main3}
\end{cor}

\noindent
\textbf{Proof}. 
Assume that $\A^H$ is not free. Then by Corollary \ref{cormain} (2), 
there is a basis $\theta_2,\ldots,\theta_\ell \in D_H(\A)$ such that 
$\theta_E^H$ can be expressed as a linear combination of $\pi(\theta_2),\ldots,\pi(\theta_\ell)$ in 
degree $\deg \theta_E^H=|\A|-|\A^H|$. Since the $b_2$-equality holds true between $\B$ and 
$H \in \B$, the same holds true in degree $|\B|-|\A^H|<|\A|-|\A^H|$, a contradiction. \owari
\medskip



\begin{rem}
The counter example to the original Orlik's conjecture in \cite{ER} says that even when the $b_2$-equation holds, 
and $\A$ free, $\A^H$ could be non-free (see Example \ref{ER1}). However, in that example, $\chi(\A^H;t)$ does not divide $\chi(\A;t)$, though 
the $b_2$-equality holds. Hence Problem \ref{Orlik} is still open.
\end{rem}


\end{document}